\documentclass[preprint,12pt]{elsarticle}

\usepackage{graphics}
\usepackage{amsmath}
\usepackage{bm}

\usepackage{miller}
\usepackage{color}
\usepackage{epstopdf}
\usepackage{amssymb}
\biboptions{comma,square,sort&compress}
\usepackage{hyperref}

\usepackage[utf8]{inputenc}
\usepackage[T2A]{fontenc}

\journal{Elsevier}

\begin{document}

\begin{frontmatter}



\title{Numerical solution of Smoluchowski coagulation equation combined with Ostwald ripening }

\author[inst1,inst2]{Robert T. Zaks}
\author[inst1,inst2]{Sergey A. Matveev}
\author[inst3]{Margarita A. Nikishina}
\author[inst3]{Dmitri V. Alexandrov}

\affiliation[inst1]{organization={Faculty of Computational Mathematics and Cybernetics, Lomonosov Moscow State University}, addressline={Leninskie Gory, bldg. 2,}, city={Moscow}, postcode={119991}, country={Russian Federation}}
\affiliation[inst2]{organization={Marchuk Institute of Numerical Mathematics of Russian Academy of Science}, addressline={Gubkin st. 8}, city={Moscow}, postcode={119333}, country={Russian Federation}}
\affiliation[inst3]{organization={Laboratory of Multi-Scale Mathematical Modeling, Department of Theoretical and Mathematical Physics, Ural Federal University}, addressline={Lenin Ave., 51}, city={Ekaterinburg}, postcode={620000}, country={Russian Federation}}

\begin{abstract}
The processes of simultaneous coagulation and Ostwald ripening of particles in the concluding stage of phase transformation are considered. We solve the integro-differential system of Smoluchowski-type kinetic and mass balance equations using a computationally efficient numerical algorithm based on low-rank matrices. We compare our numerical solutions for different initial particle-volume distributions with the universal distribution function 
for combined coagulation and Ostwald ripening. Our calculations confirm the tendency of a particulate ensemble to the universal particle-volume distribution to be approached asymptotically after a sufficiently long time, no matter what the initial particle-volume distribution might be. 
\end{abstract}




\end{frontmatter}



\section{Introduction}
\label{INTRO}
One of the fundamental problems in condensed matter physics and materials science is the theoretical description of the crystallization process in metastable and non-equilibrium media. At the same time, studying the evolution of crystal ensembles is not only of academic interest. Understanding and precise modeling of coagulation kinetics are critical for controlling microstructure and, consequently, material properties in a wide range of fields - from controlling nanoparticle characteristics in biophysics, catalysis and photonics to predicting the durability of structural alloys and creating materials for the chemical and metallurgical industries \cite{AstrucDidier,Mullin_book,HGHM,Dubrovskii,KFT,KRT,Sholokhov,EM23}.

The crystal evolution process is usually divided into three stages (initial, intermediate, and concluding \cite{PRSA}), each of which is described by specific characteristic times and different mathematical models. At the initial stage, crystallites, known as crystallization centers, begin to form in the metastable liquid (supersaturated/supercooled) and grow rapidly. This stage passes fairly quickly and is often described simultaneously with the intermediate stage, during which additional crystallization centers continue to form and crystal growth continues (see, among others, \cite{VM1990,BGY,Shift,Barlow2022,NIA,MAI_Crystals,MMAS2026,MargaritaMMAS}). As this takes place, the particles are far enough apart and do not influence each other. The Fokker-Planck-type equation and the heat/mass balance law are used to describe crystal growth at this stage. It should be noted that there are several articles that describe the transition of a crystal ensemble from the intermediate stage to the concluding one \cite{Slezov1978,PTAnuclcoars}. When crystals grow to such sizes that they begin to influence each other's growth, the crystallizing system enters the concluding stage, in which processes such as particle coagulation, described by the Smoluchowski equation, Ostwald ripening, described by the Lifshitz-Slyozov equation, and particle disintegration, described by the fragmentation equation, are possible to occur \cite{Friedlander1966,Hunt,Slezov,MakAl_OR,Crystals2022,ZhenyaMMAS}. Generally speaking, there are only a few works that describe at least two such processes simultaneously \cite{JPhys,Spicer,Foret2012,Brilliantov,JPCM_DA_2016}. This is explained by the fact that the governing equations of the concluding stage represent a complex system of integro-differential equations with moving boundaries of evolving particles/crystals (there are no general analytical/numerical methods for solving such systems).

Schumann \cite{Schumann} obtained one of the classical results on particle coagulation in the absence of other possible effects. Namely, he derived an exact analytical solution of the unsteady Smoluchowski coagulation equation with a constant coagulation kernel. The universal particle-volume distribution function that he has found is realized over long coagulation times and does not depend on the initial particle-volume distribution. He has also analytically demonstrated this conclusion for different initial distribution functions. Another fundamental result was obtained by Lifshitz and Slyozov (LS) \cite{LP,LS1961,SS1987} for the pure Ostwald ripening (coalescence) of the particles. They have derived the universal crystal-size distribution function that is established at large times in the Ostwald ripening process \cite{Relaxation}. An exact analytical solution for the particle-volume distribution function was recently derived by Makoveeva and Alexandrov \cite{HMT2023,CNSNS} when particle coagulation and Ostwald ripening occur simultaneously. This universal particle-volume distribution is also realizing for large times and has limiting transitions to the Schumann \cite{Schumann} and LS \cite{LP,LS1961,SS1987} distributions. 

The governing equations for the studied process have integro-differential form containing the integral terms related to coagulation process and partial derivatives corresponding to migration and diffusion effects. In this paper, we exploit an efficient method basing on idea of the low-rank decomposition of the coagulation kernel coefficients \cite{matveev2015fast} and classical finite-difference discretization for the partial derivatives. This approach is rather established by today and might be applied to a broad class of related models accounting the fragmentation terms \cite{zaks2025fast} and even complicated multicomponent processes \cite{matveev2025nonnegative}. 

One may also study these processes \cite{sorokin2015monte} Monte Carlo approaches that automatically construct the solution satisfying all necessary mass conservation laws. However, the stochastic approaches are less accurate than finite-difference method. This drawback limits their utilization for studies of the asymptotic properties of the particle size distributions at large times. Utilization of the low-rank deompositions for acceleration of Monte Carlo methods \cite{osinsky2024low, matveev2025low} seems to be a nice idea for future but currently there are no extensions allowing to account the migration and diffusion effects effectively within this framework.

In the present paper, we numerically solve the model of combined particle coagulation and Ostwald ripening for different initial particle-volume distributions. We show that all of them approach the universal particle-volume distribution \cite{HMT2023,CNSNS} at large times. In other words, we demonstrate that this universal distribution is independent of the initial conditions and the fact that any particulate system tends to the universal distribution over a long period of time. 

\section{Governing equations}
\label{CoagEq}
Let us consider the phase transformation of a polydisperse ensemble of particles in a metastable environment, taking into account the processes of particle coagulation and Ostwald ripening. The particle-volume distribution function $\varphi (v,t)$ satisfies the following coagulation equation complicated by Ostwald ripening \cite{Friedlander,WL}
\begin{align}
\begin{split}
\frac{\partial \varphi}{\partial t}+ \frac{\partial \mathcal{I}}{\partial v} =  
\frac{1}{2} \int\limits_0^v \mathcal{A}(w,v-w)\varphi(w,t)\varphi(v-w,t)dw \\ -\varphi(v,t)\int\limits_0^\infty \mathcal{A}(w,v)\varphi(w,t)dw ,
\label{e1}
\end{split}
\end{align}
where $v$ and $t$ are the volume and time variables, $\mathcal{I}=(dv/dt) \varphi -\mathcal{D}_v \partial\varphi /\partial v$ is the particle flow consisting of the migration $(dv/dt) \varphi$ and diffusion $\mathcal{D}_v \partial\varphi /\partial v$ flows, $\mathcal{D}_v$ is the diffusion coefficient in particle volume space. This contribution describes Ostwald ripening (coalescence) of the particles. The integral contributions in Eq. (\ref{e1}) are responsible for particle coagulation. As this takes place, the first of them defines the rate of generation of particles with volume $v$ due to collisions of lesser particles. Moreover, particles collide with each other, ensuring the rate of loss of particles with volume $v$. This process is defined by the second integral contribution in Eq. (\ref{e1}). Note that the kinetics of both processes are determined by the coagulation kernel (collision-frequency function) $\mathcal{A}(w,v)$, which depends on the physical mechanism of interaction between particles. To simplify the mathematical model, the coagulation kernel is often considered constant $\mathcal{A}(w,v)=\mathcal{A}_o={\rm const}$. The physical reason for the constant coagulation kernel is explained by the possibility of averaging of $\mathcal{A}(w,v)$ over all combinations of particle volumes $w$ and $v$ (see, among others, \cite{HMT2023,CNSNS}). Following the studies \cite{JPAMT2023,Melikhov,RandolphWhite,YuBAIv,BAM,MSMSE}, we shall assume that the diffusion coefficient is proportional to $dv/dt$, i.e. $\mathcal{D}_v = d_o(dv/dt)$ with $d_o$ being a constant defining the intensity of particle diffusion in the space of their volumes. It is important to note that the analytical definition of $\mathcal{D}_v$ is an open problem of the statistical theory of fluids \cite{LP}. Therefore, for simplicity of the analysis, we chose the above formula, which has been confirmed experimentally \cite{JPAMT2023,Melikhov,RandolphWhite,YuBAIv}. 

The mass balance condition in a saturated solution reads as follows \cite{LP}
\begin{align}
\begin{split}
c(t)-c_i +\int\limits_0^\infty w\varphi (w,t)dw =Q
\label{e2}
\end{split}
\end{align}
with $Q$ being the initial saturation, $c(t)$ and $c_i$ being the current solute concentration and the solute concentration on a planar phase transition boundary. In the following, we shall use the particle growth law as a power function of the supersaturation, i.e $dv/dt = k_g[c(t)-c_s]^\gamma$ with $k_g$ and $\gamma$ being the growth rate constants and $c_s$ being the solute concentration at saturation. This growth rate law is confirmed by experimental studies \cite{Mullin_book,Vollmer,Barlow2009,Barlow2017,JAP}.

Taking all the aforementioned aspects into account, we rewrite Eq. (\ref{e1}) as follows
\begin{align}
\begin{split}
\frac{\partial \varphi}{\partial t}+ k_g\left[ c(t)-c_s \right]^\gamma \left( \frac{\partial \varphi}{\partial v} -d_o \frac{\partial^2 \varphi}{\partial v^2} \right)  \\ 
=\frac{\mathcal{A}_o}{2} \int\limits_0^v \varphi(w,t)\varphi(v-w,t)dw -\mathcal{A}_o\varphi(v,t,)\int\limits_0^\infty \varphi(w,t)dw  .
\label{e3}
\end{split}
\end{align}

Let us especially underline that the particle-volume distribution function $\varphi (v,t)$ defines the total number $N(t)$ of particles of all volumes and the volume $V$ of new phase per unit volume of space as follows
\begin{align}
\begin{split}
N(t)=\int\limits_0^\infty \varphi (v,t)dv,\ V(t)=\int\limits_0^\infty v\varphi (v,t)dv.
\label{e4}
\end{split}
\end{align}

For the sake of simplicity, we introduce the dimensionless functions and variables in the form
\begin{align}
\begin{split}
\xi =N(0)v,\ \eta =N(0)w,\ \tau = \mathcal{A}_oN(0)t,\\ 
\Phi (\xi,\tau )=\frac{\varphi (v,t)}{N^2(0)},\ \Delta (\tau )=\frac{c(t)-c_s}{c_s},\ \kappa = \frac{k_g c_s^\gamma}{\mathcal{A}_o},\ \chi = d_o\kappa N(0).
\label{e5}
\end{split}
\end{align}
Let us underline, for definiteness, that  $\xi$ and $\eta$ represent the dimensionless particle volumes, $\tau$ represents the dimensionless time, $\Phi (\xi,\tau )$ represents the dimensionless particle-volume distribution function, and $\Delta (\tau )$ represents the dimensionless supersaturation. 


Rewriting the coagulation - Ostwald ripening model (\ref{e2})-(\ref{e4}) in dimensionless variables (\ref{e5}), we come to
\begin{align}
\begin{split}
\label{e6}
\frac{\partial \Phi}{\partial \tau}+ \Delta^\gamma (\tau )\left(\kappa\frac{\partial \Phi}{\partial \xi} -\chi \frac{\partial^2\Phi}{\partial\xi^2}\right) \\ =\frac{1}{2} \int\limits_0^\xi \Phi(\eta,\tau)\Phi(\xi -\eta,\tau)d\eta -\Phi(\xi,\tau)\int\limits_0^\infty \Phi(\eta,\tau)d\eta ,
\end{split}
\end{align}
\begin{align}
\begin{split}
V(\tau) =Q +c_i - c_s\left[\Delta (\tau) +1\right],
\label{e7}
\end{split}
\end{align}
\begin{align}
\begin{split}
n(\tau)=\frac{N(t)}{N(0)}=\int\limits_0^\infty \Phi(\xi,\tau)d\xi,\ V(\tau) = \int\limits_0^\infty \xi \Phi(\xi,\tau)d\xi.
\label{e8}
\end{split}
\end{align}

An exact analytical solution of model (\ref{e6})-(\ref{e8}) has recently been derived in Ref. \cite{CNSNS} and reads as follows (for more details, see \ref{Appendix})
\begin{align}
\begin{split}
\Delta(b) = \left\{
\begin{array}{l}
\displaystyle \left[ \Delta_o^{1-\gamma} +\frac{\kappa (1-\gamma)}{c_s}\left( \ln\left( \frac{b^2}{b_o^2} \right) +\frac{2\chi (b-b_o)}{\kappa}\right) \right]^{\frac{1}{1-\gamma}}  , \quad \gamma \neq 1\\\\
\displaystyle \Delta_o \left( \frac{b}{b_o} \right)^{\frac{2\kappa}{c_s}}\exp\left[ \frac{2\chi (b-b_o)}{c_s} \right] ,\quad \gamma =1
\end{array}
\right. ,
\label{e9}
\end{split}
\end{align}
\begin{align}
\begin{split}
n(b)=-\frac{2h(b)}{b},\ V(b)=-\frac{2h(b)}{b^2}, 
\label{e10}
\end{split}
\end{align}
\begin{align}
\begin{split}
\tau (b) = \int\limits_{b_o}^b \frac{d\beta}{h(\beta)},
\label{e11}
\end{split}
\end{align}
\begin{align}
\begin{split}
\Phi (\xi,b) = -2h(b) \exp (-b\xi ) = \frac{n^2(b)}{V(b)}\exp\left( -\frac{n(b)}{V(b)}\xi \right), 
\label{e12}
\end{split}
\end{align}
\begin{align}
\begin{split}
h(b)=-\frac{\Phi (0,0)}{2b_o^2}b^2 +b^2\int\limits_{b_o}^b \frac{\Delta^{\gamma} (\beta)(\kappa +\chi \beta)}{\beta}d \beta ,
\label{e13}
\end{split}
\end{align}
where $\Delta_o$ and $b_o$ stand for the initial values of $\Delta $ and $b$.

The analytical solution (\ref{e9})-(\ref{e13}) is found in parametric form with $b$ being the parameter. This solution defines the liquid supercooling $\Delta$ by Eq. (\ref{e9}), the dimensionless total number $n$ of particles and the volume $V$ of the condensed phase by Eqs. (\ref{e10}) and time $\tau$ by Eq. (\ref{e11}) as functions of the parameter $b$. In other words, knowing the dependence $\tau (b)$, one can construct the inverse function $b(\tau)$ and establish the dependence of all unknown functions on the dimensionless time $\tau$. As this takes place, the particle-volume distribution function $\Phi$ is defined by Eq. (\ref{e12}). An important point is that this function defines the universal distribution by particle volumes, which is established at large times (see, for more details, discussions in Refs. \cite{HMT2023,CNSNS,Schumann}). 

Note that the distribution function (\ref{e12}) satisfies the boundary conditions
\begin{equation}\label{eq:boundary}
    \Phi\xrightarrow[\xi\to\infty]{}0,\,
    \dfrac{\partial\Phi}{\partial\xi}\xrightarrow[\xi\to\infty]{}0,
\end{equation}
which show that $\Phi$ gradually decays at high particle volumes.

\section{Numerical method and experiments}
\label{NummericalMethod}

\subsection{Discretization of the integral operators and time-integration}

Let us suppose that $\operatorname{supp}\mathcal{A(\xi, \eta)} 
\subseteq[0;H]^{2}$ and
to find the distribution function $\Phi$ at time $\tau=T$, we use the uniform grids
$\{\xi_{i}=ih_{\xi},i\in\overline{0,M}\}$,
$\{\tau_{j}=jh_{\tau},n\in\overline{0,M_{\tau}}\}$
within the intervals $[0;H]$ and $[0;T]$.

First of all, we utilize the simplest explicit Euler time-integration method and introduce the integral orerators $L_1(\Phi)(\xi)$ and $L_2(\Phi)(\xi)$:
\begin{align}
\begin{split}
\label{eq:Euler_step}
& \frac{\partial \Phi}{\partial \tau} \approx \frac{\Phi(\tau_{j+1}, \xi) - \Phi(\tau_j, \xi)}{\tau_j} = \\ =
& \underbrace{\frac{1}{2} \int\limits_0^\xi \mathcal{A}(\eta, \xi-\eta) \Phi(\eta,\tau_j)\Phi(\xi -\eta,\tau_j)d\eta}_{L_1(\Phi)(\xi)} -\Phi(\xi,\tau_j)\underbrace{\int\limits_0^\infty \mathcal{A}(\xi, \eta)\Phi(\eta,\tau_j)d\eta}_{L_2(\Phi)(\xi)}
- \\
& - \underbrace{\Delta^\gamma (\tau_j )\left(\kappa\frac{\partial \Phi}{\partial \xi} -\chi \frac{\partial^2\Phi}{\partial\xi^2}\right)}_{\text{partial derivatives at moment $\tau_j$}},
\end{split}
\end{align}
Fixing the time moment $\tau_j$ (below, we skip this variable in equation \eqref{eq:OneDimensionalTrapezoids} for simpler notations), we further utilize the uniform grid over $\xi$ axis and approximate these integrals via the trapezoid rule:
\begin{small}
\begin{equation}\label{eq:OneDimensionalTrapezoids}
\begin{matrix}
L_1 \Phi (\xi_i) = \frac{h_{\xi}}{2} \sum\limits_{i_1=1}^i  \left( \mathcal{A}( \xi_{i-i_1} , \xi_{i_1}) \Phi(\xi_{i_1}) \Phi(\xi_{i-i_1})  + \mathcal{A}(\xi_{i-i_1+1} , \xi_{i_1 -1}) \Phi(\xi_{i_1 - 1} ) \Phi(\xi_{i-i_1+1})  \right) , \\ i = 1, 2 , \ldots, M, ~ L_1 \Phi (0) = 0;  \\ 
L_{2} \Phi (\xi_i)= \frac{h_{\xi}}{2}  \mathcal{A}(\xi_i, \xi_0) \Phi(\xi_0)  + \frac{h_{\xi}}{2}  \mathcal{A}(\xi_i, \xi_M) \Phi(\xi_M) + h_{\xi} \sum\limits_{i_1=1}^{M-1}  \mathcal{A}(\xi_i, \xi_{i_1}) \Phi(\xi_{i_1}) , \\ i = 0, 1, 2 \ldots,  M.
\end{matrix}
\end{equation}
\end{small}
If the coagulation kernel admits the low-rank representation
\begin{equation}\label{eq:SkelDecomp}
\mathcal{A}(\xi, \eta) = \sum_{\alpha=1}^{R} u_{\alpha}(\xi) \cdot v_{\alpha}(\eta)
\end{equation}
with $R \ll M$, the integrals from equation \eqref{eq:OneDimensionalTrapezoids} can be evaluated within the $O(R M \log M)$ operations (see \cite{matveev2015fast} for detail) via combination of matrix-based operations for $L_2(\Phi)(\xi)$ and FFT for $L_1(\Phi)(\xi)$ as lower-triangular convolution. In order to implement the final algorithm for the time-step, we need to discuss the discretization of the partial derivatives in \eqref{eq:Euler_step}. In the current paper, we study the problem with the constant kernel $\mathcal{A}(\xi, \eta) \equiv A_0=  {\rm const}$ having rank 1 and meaning that this methodology allows to make dramatically efficient calculations.

The class of kernels having a good form in this format is rather broad. For example, the famous diffusion kernel \cite{Friedlander,WL}
$$
\mathcal{A}(\xi, \eta) = \left( \xi^{1/3} + \eta^{1/3} \right) \left(\xi^{-1/3} + \eta^{-1/3}\right) = 
\left( \frac{\xi}{\eta}\right)^{1/3} + ~\left( \frac{\eta}{\xi} \right)^{1/3} +~ 2
$$
has an exact representation in such a form with number of terms equal to 3. Another interesting example is the free-molecular kernel (that is also known as ballistic) \cite{Friedlander,WL}
$$
\mathcal{A}(\xi, \eta) = \left( \xi^{1/3} + \eta^{1/3} \right)^2 \sqrt{\frac{1}{\xi} + \frac{1}{\eta}},
$$
that does not have an exact representation in the low-rank form \eqref{eq:SkelDecomp} but can be approximated with $O(\log \frac{1}{\varepsilon}) $ terms, where $\varepsilon$ is the accuracy of the corresponding approximation (see the proof in \cite{matveev2015fast}). 

\subsection{Discretization of the partial derivatives}

In addition, we introduce an additional absorption layer (see e.g. \cite{berenger1994perfectly, zagidullin2017efficient})
$\{\xi_{i}=ih_{\xi},i\in\overline{M+1,M_{\xi}}\}$
to satisfy the boundary conditions \eqref{eq:boundary}. Here, after each step, the calculated values $\Phi(\xi_{i},\tau_{n}),i\in\overline{M+1,M_{\xi}}$ are multiplied by $\exp\left(-d(i-M)h_{\xi}\right)$ with $d$ being a parameter that is chosen to set $\Phi$ and
$\dfrac{\partial\Phi}{\partial\xi}$ zero at the absorption layer boundary. Note that the choice of exponential dependence is not accidental: stationary
solutions of equation \eqref{e6} (at least for $\chi=0$, $\kappa=0$)
have precisely this asymptotic dependence for $\xi\to\infty$.

We approximate the derivatives using the following relations for numerical differentiation:
\begin{gather*}
    \dfrac{\partial\Phi}{\partial\tau}(\xi_{i},\tau_{j+1})=
    \dfrac{\Phi_{i}^{j+1}-\Phi_{i}^{j}}{h_{\tau}}+O(h_{\tau}),\\
    \dfrac{\partial\Phi}{\partial\xi}(\xi_{i},\tau_{j+1})=
    \begin{cases}
        \dfrac{-3\Phi_{0}^{j}+4\Phi_{1}^{j}-\Phi_{2}^{j}}{h_{\xi}},&i=0\\
        \dfrac{\Phi_{i+1}^{j}-\Phi_{i-1}^{j}}{2h_{\xi}},
        &i\in\overline{1,M_{\xi}-1}\\
        \dfrac{
            -\Phi_{M_{\xi}-2}^{j}+4\Phi_{M_{\xi}-1}^{j}-
            3\Phi_{M_{\xi}}^{j}
        }{h_{\xi}},&i=M_{\xi}
    \end{cases}+O(h_{\xi}^{2}),\\
    \dfrac{\partial^{2}\Phi}{\partial\xi^{2}}(\xi_{i},\tau_{j+1})=
    \begin{cases}
        \dfrac{2\Phi_{0}^{j}-5\Phi_{1}^{j}+4\Phi_{2}^{j}-\Phi_{3}^{j}}{
            h_{\xi}^{2}
        },&i=0\\
        \dfrac{\Phi_{i+1}^{j}-2\Phi_{i}^{j}-\Phi_{i-1}^{j}}{h_{\xi}^{2}},
        &i\in\overline{1,M_{\xi}-1}\\
        \dfrac{
            -\Phi_{M_{\xi}-3}^{j}+4\Phi_{M_{\xi}-2}^{j}-
            5\Phi_{M_{\xi}-1}^{j}+2\Phi_{M_{\xi}}^{j}
        }{
            h_{\xi}^{2}
        },&i=M_{\xi}
    \end{cases}+O(h_{\xi}^{2})
\end{gather*}
One may also exploit the finite-volume framework \cite{singh2016volume} as alternative methodology better matching the conservation laws. However, in this work we use the finite-difference approach allowing to obtain the second order of accuracy as well as the low-rankness of the coagulation kernel for fast computing. The finite-volume framework would cause change in calculations of the coagulation integrals making them slower (see e.g. \cite{osinsky2020low} containing detailed benchmarks) that is unpleasant for the studies of numerical solutions for long times

Finally, to approximate the integral operator in the right-hand side of equation
\eqref{e6}, as we discussed above, we use the specific implementation of our numerical method from \cite{zaks2025fast}. More specifically, at first, the approximation of the integral operator in the right-hand side of equation \eqref{e6} is applied, then the values of $\Phi_{i}^{j+1}$ are calculated, and then the absorbing layer is applied. This numerical scheme is implemented as a program using the \texttt{SmoluchowskiSolver} library and is openly available by link \url{https://smoluch.rozax.net/download/SmoluchowskiSolver\_partial.tgz}\footnote{\label{downloadTGZ}To unzip the archive, use the command:\\
\texttt{tar czf SmoluchowskiSolver\_partial.tgz},
then, with \texttt{g++} and \texttt{gnu make} installed,
you need to write \texttt{make} and go to the directory \texttt{./examples} and
run the program \texttt{./partial}.}.

\subsection{Verification of numerical method}

We verify our numerical method comparing the analytical and numerical solutions of equations \eqref{e6}--\eqref{e8} in the case of $\mathcal{A}=\mathcal{A}_o = {\rm const}$ and $b$ changing from $b=b_{o}=1$ to $b=b_{1}=0$. Let us note several features of the solving procedures: (i) to calculate the distribution function (\ref{e12}), it is necessary to find numerically the values of two integrals, and (ii) to calculate $\Phi$ at time $T$, it is necessary to select the values of the parameter $b$ for which $T=\tau(b)$.

\begin{figure}[!t]
    \centering
    \begin{minipage}[b]{0.495\textwidth}
        \includegraphics[width=\textwidth]{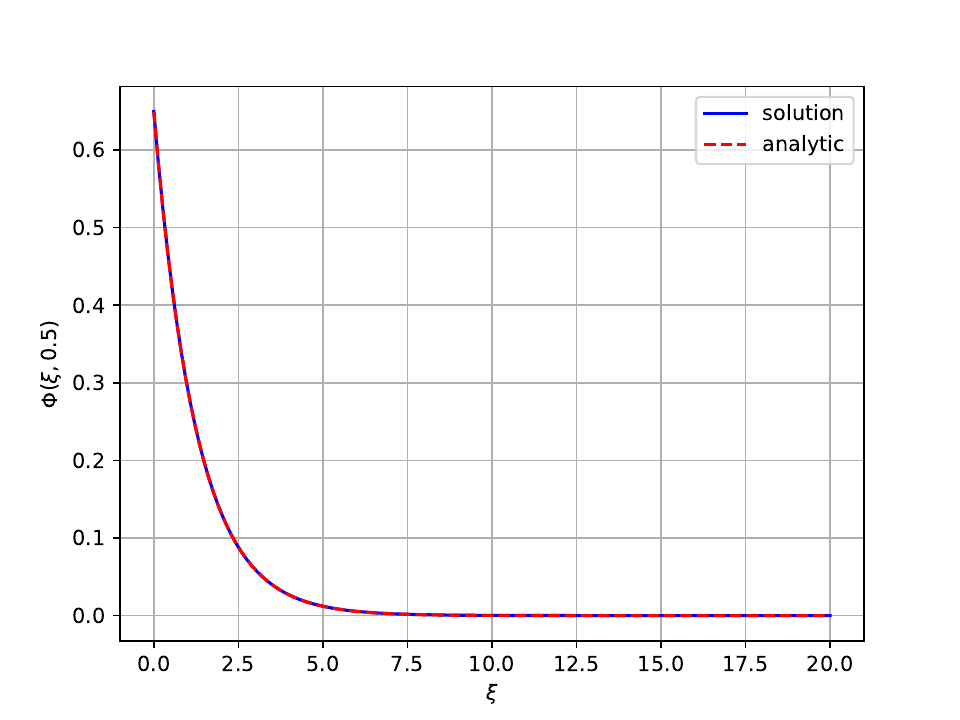}
    \end{minipage}
    \begin{minipage}[b]{0.495\textwidth}
        \includegraphics[width=\textwidth]{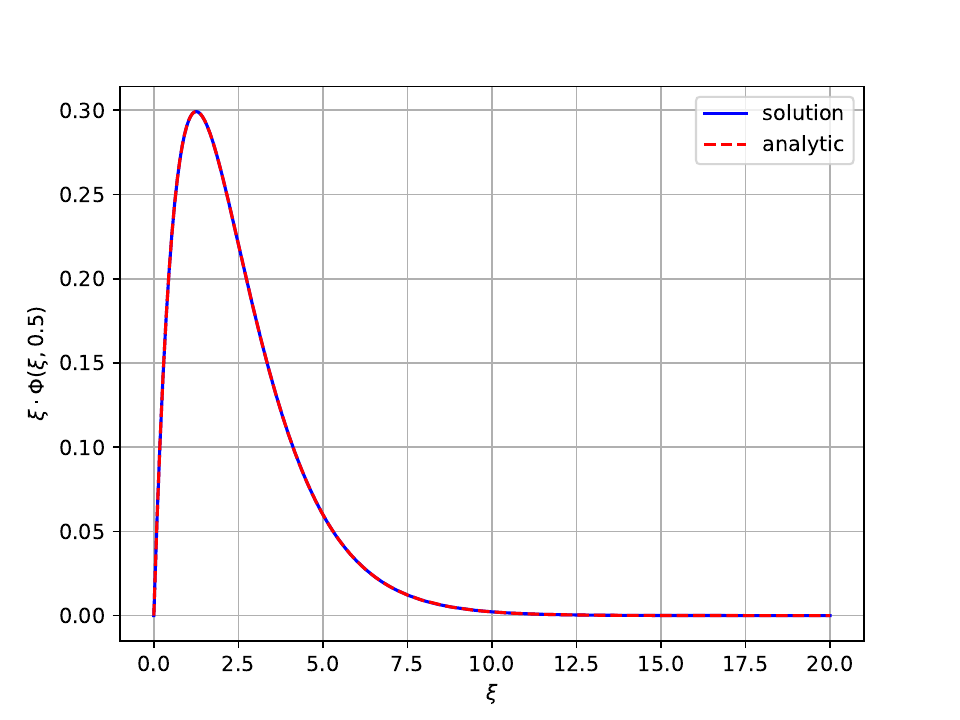}
    \end{minipage}
    \caption{\centering
Comparison of analytical and numerical solutions for $\gamma=1$, $\kappa=0.2$,
$\chi=10^{-2}$, $\Delta_{o}=0.2$, $\Phi(0,0)=1$, $b_{o}=1$ and $c_{s}=10$ mol m$^{-3}$ \cite{CNSNS}. For numerical calculations the following parameters have also been chosen:
${H=20}$, $M=4000$, $M_{\xi}=M$, $T=0.5$, $M_{\tau}=20000$ and $d=5$.
    }
    \label{fig:gamma1}
\end{figure}

\begin{figure}[!ht]
    \centering
    \begin{minipage}[b]{0.495\textwidth}
        \includegraphics[width=\textwidth]{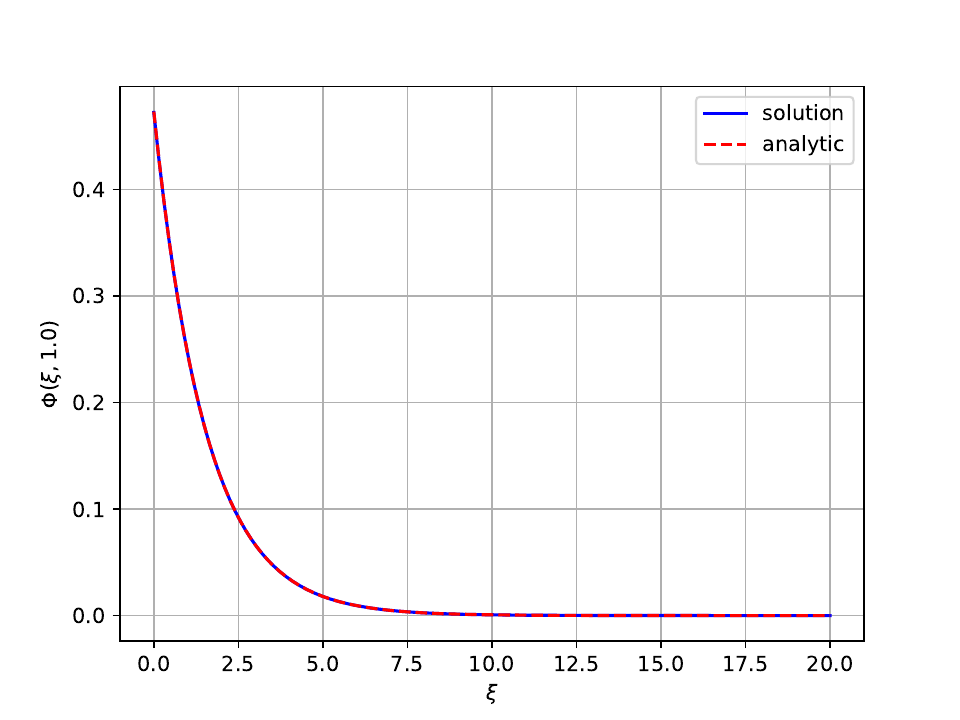}
    \end{minipage}
    \begin{minipage}[b]{0.495\textwidth}
        \includegraphics[width=\textwidth]{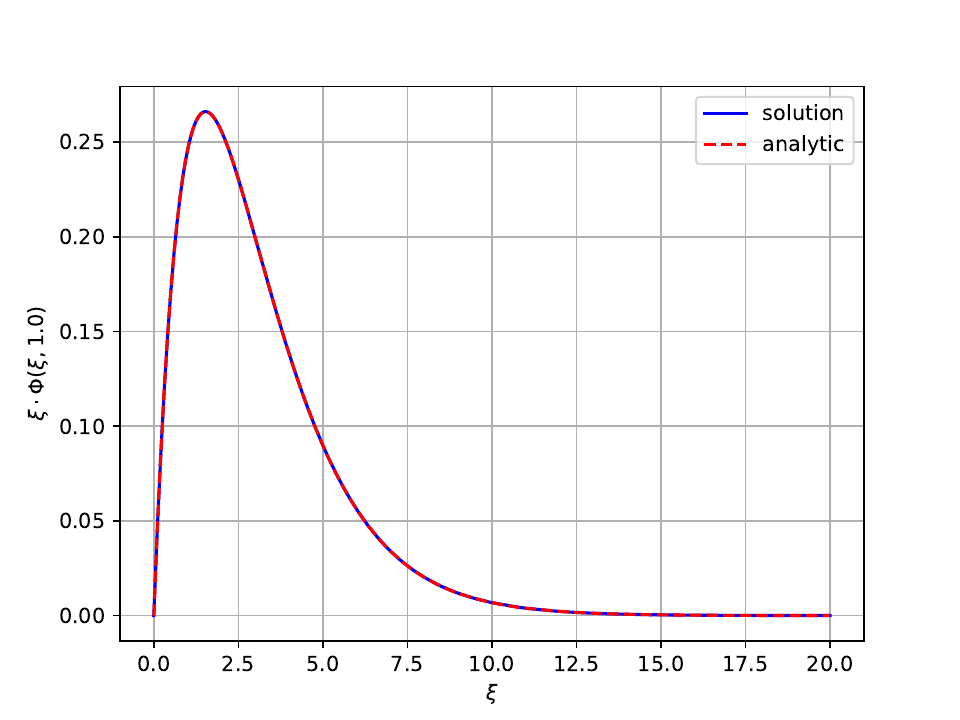}
    \end{minipage}
    \caption{\centering
        Comparison of analytical and numerical solutions for $\gamma=0.5$, $\kappa=0.2$,
$\chi=10^{-1}$, $\Delta_{o}=0.2$, $\Phi(0,0)=1$, $b_{o}=1$ and $c_{s}=10$ mol m$^{-3}$ \cite{CNSNS}. For numerical calculations the following parameters have also been chosen:
$H=20$, $M=4000$, $M_{\xi}=M$, $T=1.0$, $M_{\tau}=40000$ and $d=5$.
    }
    \label{fig:gamma05}
\end{figure}

Figures \ref{fig:gamma1} and \ref{fig:gamma05} compare the analytical particle-radius distribution function $\Phi$ and the product $\xi\Phi$ calculated according to expression (\ref{e12}) with the corresponding functions found numerically. As this takes place, the integrals of the analytical solution have been evaluated using the trapezoidal rule with a uniform grid of $40,000$ nodes. As is easily seen, both solutions coincide, i.e. the numerical scheme is verified by means of the aforementioned exact analytical solution.


\subsection{Particle-volume distribution at large time}
\label{Res}

Let us consider the numerical solution of the problem for various initial distribution functions and compare it with the universal distribution function (\ref{e12}) derived analytically. For definiteness, let us assume that $\Phi(0,0)=1$, $b_{o}=1$ and consider four different initial particle-volume distributions listed in Table \ref{tab1}.

\begin{table}[h!]
    \centering
\begin{tabular}{c|c}
        & \\
        Initial particle-volume distribution & Designation \\
        & \\
        \hline
        & \\
        $\Phi(\xi,0)=\Phi(0,0)\exp(-b_{o}\xi)$ &  \texttt{exp}  \\
        & \\
        $\Phi(\xi,0)=C_{1}\exp(-1.6\cdot b_{o}\xi)$ & \texttt{pert exp}  \\
        & \\
        $\Phi(\xi,0)=C_{2}\exp(-\xi^{2})$ & \texttt{gaus}          \\
        & \\
         $\Phi(\xi,0)=C_{3}
        \exp{\left(-\dfrac{(\xi-2)^{2}}{4}\right)}$   &   \texttt{gaus2}
\end{tabular}
    \caption{ Initial particle-volume distributions used in calculations. Here, the constants $C_{1}$, $C_{2}$, and $C_{3}$ are chosen so that the initial values of $V(0)$ for all cases under question are coincided. }
    \label{tab1}
\end{table}

\begin{figure}[!ht]
    \centering
    \begin{minipage}[b]{0.495\textwidth}
        \includegraphics[width=\textwidth]{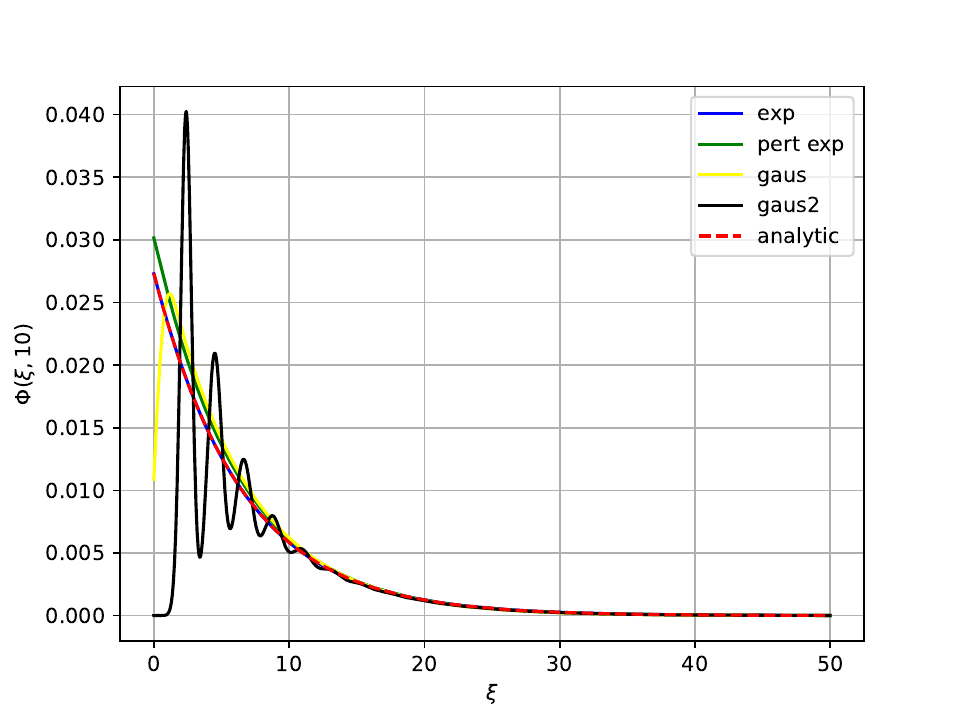}
    \end{minipage}
    \begin{minipage}[b]{0.495\textwidth}
        \includegraphics[width=\textwidth]{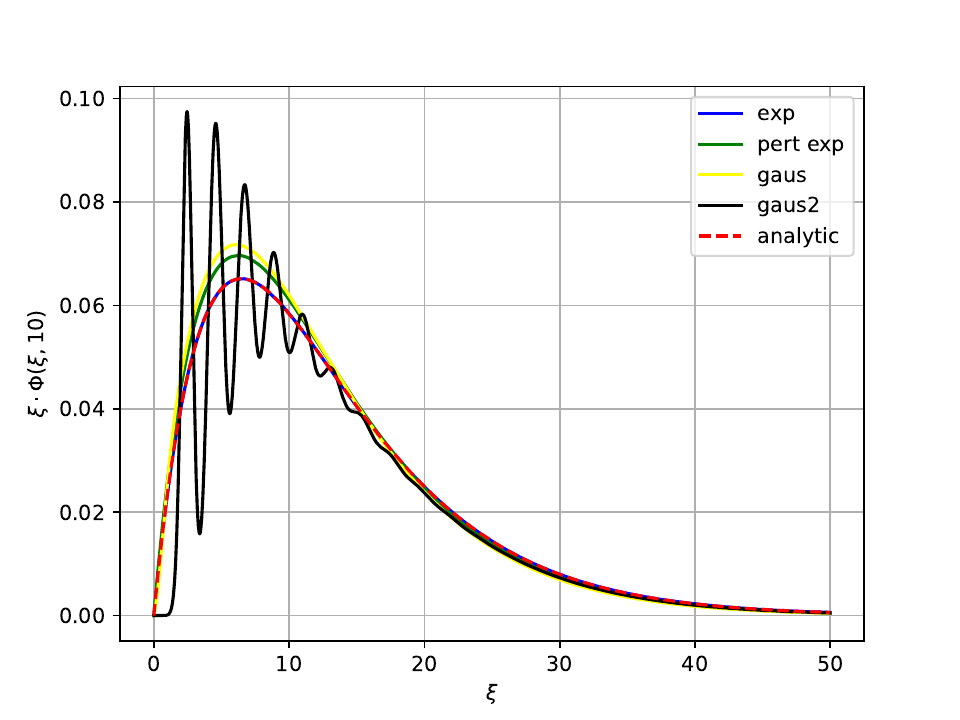}
    \end{minipage}
    \caption{\centering
        Numerical solutions for $\Phi$ and $\xi\Phi$ with various initial particle-volume distributions at $T=10$ versus the universal distribution function (\ref{e12}).
    }
    \label{fig:T10}
\end{figure}

\begin{figure}[!ht]
    \centering
    \begin{minipage}[b]{0.495\textwidth}
        \includegraphics[width=\textwidth]{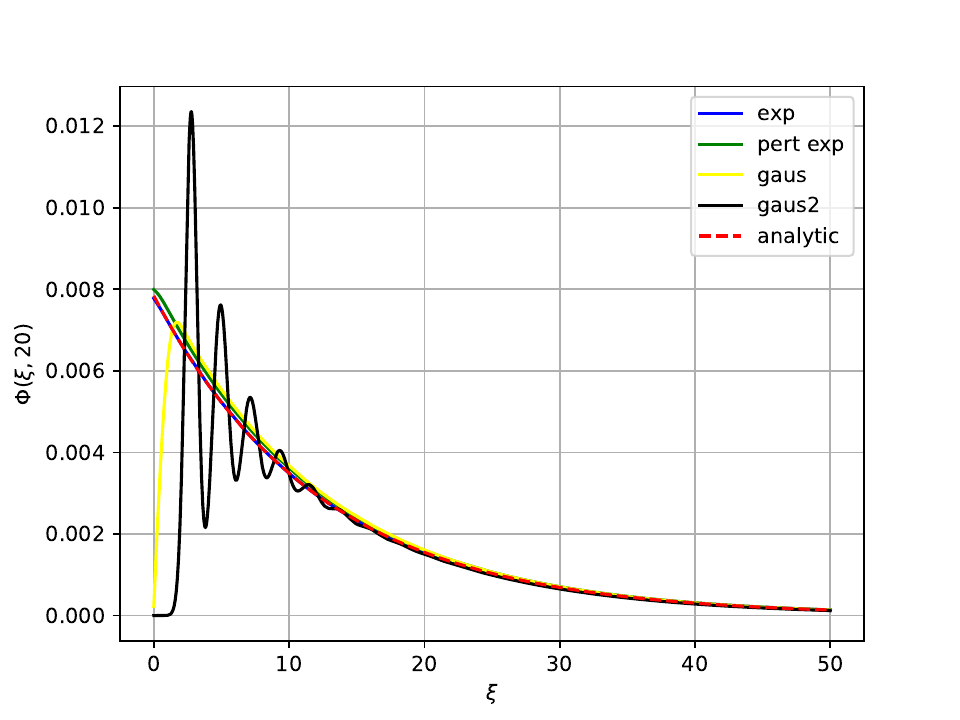}
    \end{minipage}
    \begin{minipage}[b]{0.495\textwidth}
        \includegraphics[width=\textwidth]{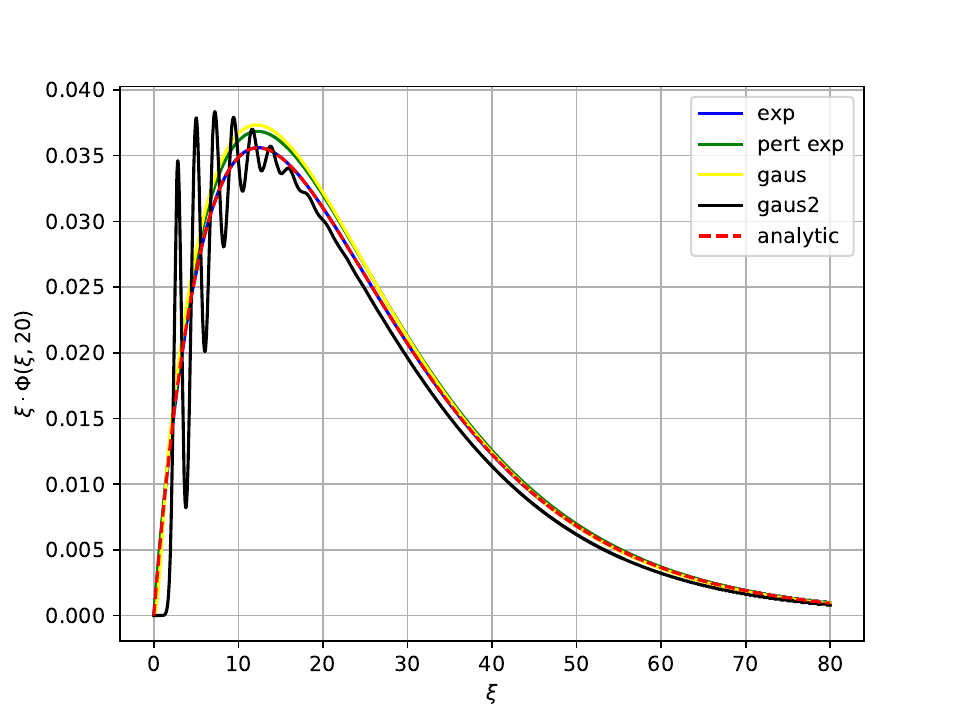}
    \end{minipage}
    \caption{\centering
        Numerical solutions for $\Phi$ and $\xi\Phi$ with various initial particle-volume distributions at $T=20$ versus the universal distribution function (\ref{e12}).
    }
    \label{fig:T20}
\end{figure}

\begin{figure}[!ht]
    \centering
    \begin{minipage}[b]{0.495\textwidth}
        \includegraphics[width=\textwidth]{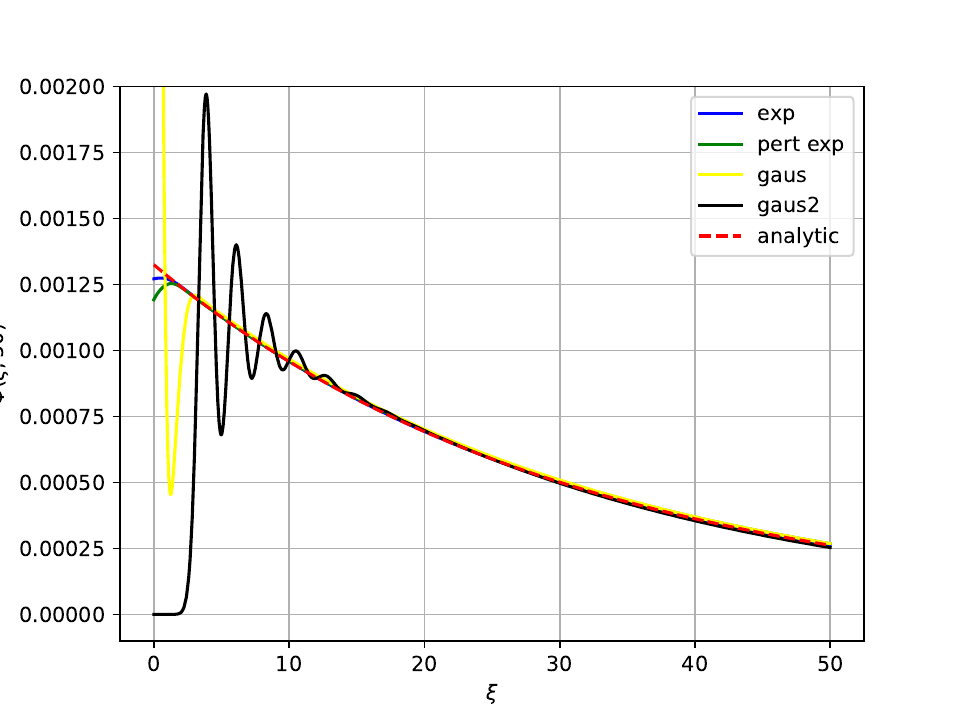}
    \end{minipage}
    \begin{minipage}[b]{0.495\textwidth}
        \includegraphics[width=\textwidth]{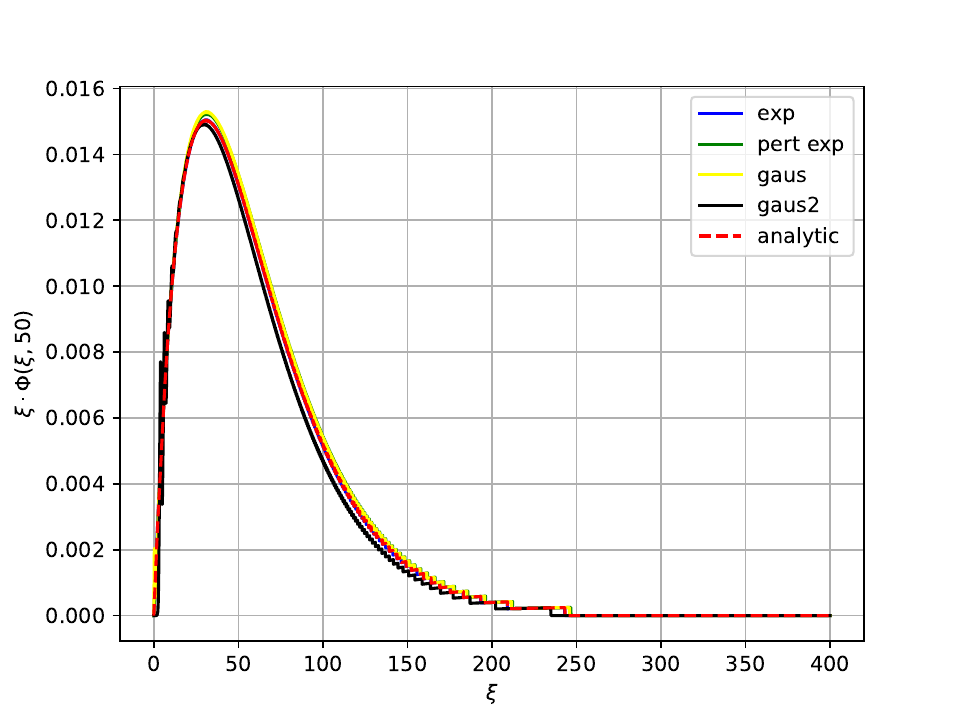}
    \end{minipage}
    \caption{\centering
        Numerical solutions for $\Phi$ and $\xi\Phi$ with various initial particle-volume distributions at $T=50$ versus the universal distribution function (\ref{e12}).
    }
    \label{fig:T50}
\end{figure}

Our numerical solutions in comparison with the analytical particle-volume distribution (\ref{e12}) are illustrated in Figs. \ref{fig:T10}, \ref{fig:T20} and \ref{fig:T50} for the following parameters $\gamma=1$, $\kappa=0.2$,
$\chi=10^{-2}$, $\Delta_{o}=0.2$, $\Phi(0,0)=1$, $b_{o}=1$, $c_{s}=10$ mol m$^{-3}$, $H=400$, $M=100\cdot H$, $M_{\xi}=M$, $T=10$,
$M_{\tau}=T\left(\dfrac{M}{H}\right)^{2}$.

It can be seen that the instability of the proposed numerical scheme already manifests itself for $T=50$ (for small values of $\xi$
there is a noticeable difference between the numerical solution \texttt{exp} and the analytical solution
\texttt{analytic}). As $T$ increases, this discrepancy will accumulate; however, to simplify implementation, we can consider precisely this numerical scheme, limiting ourselves to small values of $T$.

\begin{figure}[!ht]
    \centering
    \begin{minipage}[b]{0.495\textwidth}
        \includegraphics[width=\textwidth]{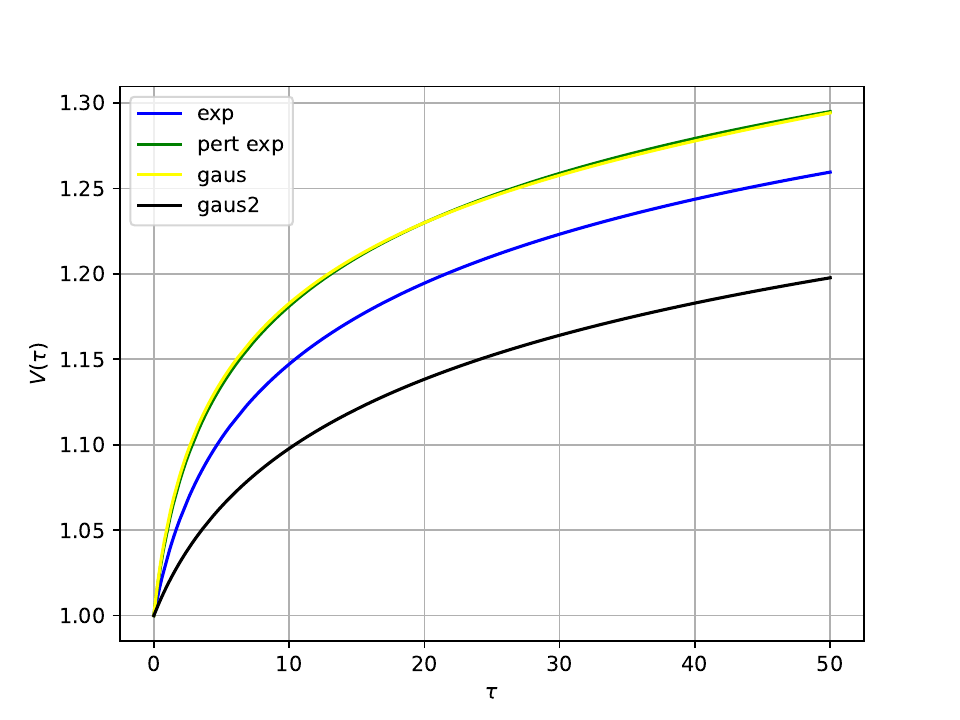}
    \end{minipage}
    \begin{minipage}[b]{0.495\textwidth}
        \includegraphics[width=\textwidth]{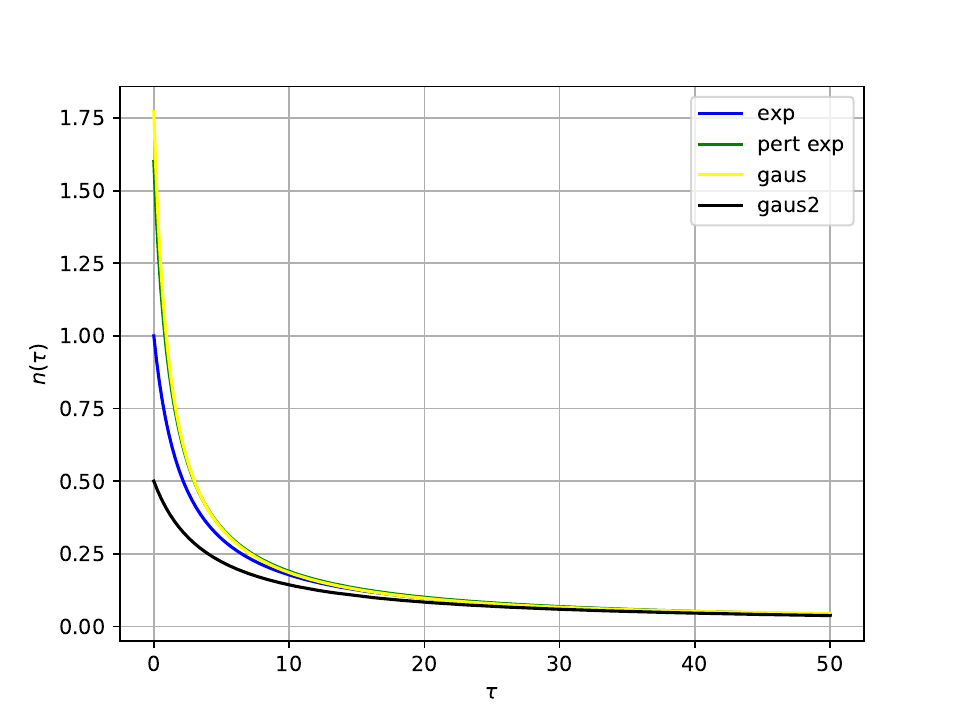}
    \end{minipage}
    \caption{\centering
        The volume $V$ of the new phase and reduced total number $n$ of particles of all volumes as functions of dimensionless time $\tau$ for $H=400$.
    }
    \label{fig:md}
\end{figure}

Figs. \ref{fig:T10}-\ref{fig:T50} show that all different initial particle-volume distributions lead to the universal distribution function (\ref{e12}) for $\xi\gtrsim 10$. {\it This confirms a tendency to the same universal particle-volume distribution (\ref{e12}) to be approached asymptotically after a sufficiently long time, no matter what the initial particle-volume distribution might be}. A strict proof of this conclusion appears to be very difficult, and we must content ourselves with the demonstration that at least several forms of particle-volume initial distributions lead to the universal distribution function (\ref{e12}), taking into account simultaneous occurrence of coagulation and Ostwald ripening mechanisms. An important point is that the universal particle-volume distribution was found by Schumann \cite{Schumann} in the case of pure coagulation, and by Lifshitz and Slyozov in the case of pure Ostwald ripening (coalescence) \cite{LP,LS1961,SS1987}. Therefore, their combined effect (coagulation and Ostwald ripening), as our calculations show, should have the same tendency to the universal particle-volume distribution (\ref{e12}) derived in Refs. \cite{HMT2023,CNSNS}. Fig. \ref{fig:md} shows the same large-time behavior for the volume $V$ of the new phase and reduced total number $n$ of particles of all volumes: 
$\displaystyle V(\tau)\approx\int\limits_{0}^{H}{\xi \Phi(\xi,\tau)d\xi}$,
$\displaystyle n(\tau)\approx\int\limits_{0}^{H}{\Phi(\xi,\tau)d\xi} $.

\section{Concluding remarks}
\label{Conc}

In summary, we briefly outline the basic conclusions that follow from our study. We consider a concluding stage of the phase transformation process when the coagulation of particles is complicated by their Ostwald ripening. The mathematical model consists of the integro-differential Smoluchowski-type kinetic equation and the integral mass balance equation supplemented with corresponding initial and boundary conditions. The problem is solved numerically using the openly available \texttt{SmoluchowskiSolver} (see the link above). 

The proposed numerical scheme for the coagulation
equation combined with Ostwald ripening in the case of a constant coagulation kernel gives an approximation that agrees with the analytical solution recently found in Ref. \cite{CNSNS} and describing the universal particle-volume distribution at large times. Our numerical calculations demonstrate that different initial particle-volume distributions lead a particulate ensemble to the same universal distribution function (\ref{e12}) for crystal volumes $\xi\gtrsim 10$. {\it This fact confirms the tendency of a particulate ensemble to the universal particle-volume distribution (\ref{e12}) to be approached asymptotically after a sufficiently long time, no matter what the initial particle-volume distribution might be}. So, for example, considering pure particle coagulation, the distribution function (\ref{e12}) tends to the Schumann distribution \cite{Schumann}, while considering pure Ostwald ripening, the distribution function (\ref{e12}) tends to the Lifshitz-Slyozov distribution \cite{LP,LS1961,SS1987} (fore more details, see also \cite{HMT2023,CNSNS}).

The present theory can be extended in future studies to account for the following tasks: (i) To consider different mechanisms of particle coagulation, one can use averaged coagulation kernels \cite{HMT2023} and the present numerical scheme. (ii) To consider particle fragmentation together with their coagulation and Ostwald ripening, one can use the same approach consisting of the theory of Refs. \cite{HMT2023,CNSNS} and the present numerical technique. (iii) Numerical schemes for approximating derivatives with greater stability can be considered. (iv) Implementation and cross-validation of the extenstions of the efficient stochastic or finite-volume approaches accounting the migration and diffusion terms.

\section*{CRediT authorship contribution statement}

Robert T. Zaks:
Conceptualization, Data curation, Formal analysis, Investigation, Funding acquisition, Methodology, Project administration, Resources, Software, Supervision, Validation, Visualization, Writing - original draft, Writing - review and editing.

Sergey A. Matveev:
Conceptualization, Data curation, Formal analysis, Investigation, Funding acquisition, Methodology, Project administration, Resources, Software, Supervision, Validation, Visualization, Writing - original draft, Writing - review and editing.

Margarita A. Nikishina:
Conceptualization, Data curation, Formal analysis, Investigation, Funding acquisition, Methodology, Project administration, Resources, Software, Supervision, Validation, Visualization, Writing - original draft, Writing - review and editing.

Dmitri V. Alexandrov: Conceptualization, Data curation, Formal analysis, Investigation, Funding acquisition, Methodology, Project administration, Resources, Software, Supervision, Validation, Visualization, Writing - original draft, Writing - review and editing.

\section*{Declaration of Competing Interest}
The authors declare no conflict of interest.

\section*{Acknowledgments}
This study was financially supported by the Ministry of Science and High Education of the Russian Federation (Ural Federal University Program of Development within the Priority-2030 Program).


\appendix

\section{Analytical solution to Smoluchowski's coagulation equation combined with Ostwald ripening}
\label{Appendix}
Let us demonstrate how to derive an exact analytical solution of dimensionless equations (\ref{e6})-(\ref{e8}) that describe synchronous coagulation and Ostwald ripening.  

First, we multiply Eq. (\ref{e6}) by $d\xi$ and $\xi d\xi$ and integrate the result from $\xi =0$ to $\xi\to\infty$ with allowance for boundary conditions (\ref{eq:boundary}). As a result, we come to the equations for $n(\tau)$ and $V(\tau)$ (for a more detailed explanation, we refer the reader to the classic monograph by Williams and Loyalka \cite{WL}, Section 3.2):
\begin{align}
\begin{split}
\frac{dn}{d\tau} - \Delta^\gamma (\tau )\left(\kappa \Phi (0,\tau) -\chi  \frac{\partial \Phi (0,\tau)}{\partial\xi}\right) +\frac{n^2(\tau)}{2}=0,
\label{8}
\end{split}
\end{align}
\begin{align}
\begin{split}
\frac{dV}{d\tau} - \Delta^\gamma (\tau) \left( \kappa n(\tau) +\chi \Phi (0,\tau)\right)  =0.
\label{9}
\end{split}
\end{align}

Second, substituting $V(\tau)$ from Eq. (\ref{e7}) to Eq. (\ref{9}), we obtain an equation containing the total number $n(\tau)$ of particles, liquid supersaturation $\Delta (\tau)$ and particle-volume distribution function $\Phi (0,\tau)$ as follows
\begin{align}
\begin{split}
n(\tau)=-\frac{1}{\kappa}\left(\frac{c_s}{\Delta^{\gamma} (\tau)}\frac{d\Delta}{d\tau} +\chi \Phi (0,\tau )\right) .
\label{10}
\end{split}
\end{align}

Third, by analogy with the classical paper by Schumann \cite{Schumann}, we shall seek the particle-volume distribution in the following exponential form (such a solution is established for large times when an ensemble of particles has forgotten its initial state):  
\begin{align}
\begin{split}
\Phi (\xi,\tau)= \Phi (0,\tau)\exp\left[ -b(\tau)\xi \right] ,
\label{11}
\end{split}
\end{align}
where $\Phi (0,\tau)$ and $b(\tau)>0$ are found below. 

By substituting (\ref{11}) into (\ref{e6}), equating the terms with the same power functions of $\xi$, we get two equations
\begin{align}
\begin{split}
\Phi (0,\tau) =-2\frac{db}{d\tau},
\label{12}
\end{split}
\end{align}
\begin{align}
\begin{split}
\frac{d\Phi(0,\tau)}{d\tau}-\Delta^\gamma (\tau) b(\tau)\Phi (0,\tau)\left( \kappa +\chi 
 b(\tau) \right) +\frac{\Phi^2(0,\tau)}{b(\tau )}=0.
\label{13}
\end{split}
\end{align}
In addition, combining (\ref{11}) and (\ref{e8}), we have
\begin{align}
\begin{split}
n(\tau)=\frac{\Phi (0,\tau)}{b(\tau)},\ V(\tau)=\frac{\Phi (0,\tau)}{b^2(\tau)}.
\label{14}
\end{split}
\end{align}

Fourth, eliminating $\Phi (0,\tau)$ and $n(\tau)$ from Eqs. (\ref{10}), (\ref{12}) and (\ref{14}), we arrive at an equation connecting the functions $b(\tau)$ and $\Delta (\tau)$ as follows
\begin{align}
\begin{split}
2\left(\frac{1}{b(\tau)}+\frac{\chi}{\kappa}\right)\frac{db}{d\tau}=\frac{c_s}{\kappa \Delta^{\gamma}}\frac{d\Delta}{d\tau}.
\label{15}
\end{split}
\end{align}

Fifth, this equation can be easily integrated, after multiplying by $d\tau /db$. After integration, we come to an explicit dependence of the liquid supersaturation $\Delta $ on the parameter $b$ given by Eq. (\ref{e9}).

Sixth, combining Eqs. (\ref{12}) and (\ref{13}), we get the second equation containing $\Delta$ and $b$ as follows
\begin{align}
\begin{split}
\frac{d^2b}{d\tau^2}=\frac{db}{d\tau}\left[ \Delta^\gamma (\tau) b(\tau) \left( \kappa +\chi b(\tau )\right) + \frac{2}{b}\frac{db}{d\tau} \right] .
\label{17}
\end{split}
\end{align}

Now, introducing the new function $h(b) =db/d\tau \neq 0$, we have 
\begin{align}
\begin{split}
\frac{dh}{db}= \Delta^\gamma b \left( \kappa +\chi b\right) + \frac{2h(b)}{b}  .
\label{18}
\end{split}
\end{align}
Integrating this equation, we come to the explicit dependence of the function $h$ on the parameter $b$ given by Eq. (\ref{e13}). 

Seventh, integrating equation $d\tau /db =1/h(b)$, we obtain the dependence of time $\tau$ on the parameter $b$ given by Eq. (\ref{e11}).

Eighth, expressions (\ref{11}) and (\ref{14}) lead to (\ref{e10}) and (\ref{e12}), defining the reduced total number $n$ of particles, volume $V$ of the new phase and particle-volume distribution function $\Phi$ as functions of the parameter $b$.

Thus, we arrive at an exact solution (\ref{e9})-(\ref{e13}) in parametric form with $b$ being the decision variable (parameter).



\section*{Data availability}
No data was used for the research described in the article.

 \bibliographystyle{elsarticle-num} 
 \bibliography{main}





\end{document}